# EXAMINING THE POSITION OF BUILDING INFORMATION MODELING (BIM) TECHNOLOGY IN DIFFERENT DIMENSIONS OF BUILDING SMARTNESS

Paniz Pouryaghoubi[1], Arefeh Mohammadi[2*]

[1]Department of Civil Engineering, Ferdowsi University of Mashhad, Mashhad, Iran
[2]Department of Civil Engineering, Ferdowsi University of Mashhad, Mashhad, Iran

**Abstract.** Due to the fact that the construction industry is one of the high-risk industries, it is inevitable to pay sufficient attention to the safety category in this industry. Currently, due to the development of modern technologies, the traditional methods used in safety management have given way to the mentioned technologies. Therefore, the use of modern technologies such as BIM can be a practical step in reducing human and financial losses while increasing productivity and improving construction quality. The purpose of this research is to investigate the role and place of using building information modeling in different sectors, such as improving smartness, safety, and its application in reducing energy consumption. Considering that a significant share of the capital of different countries is spent in the construction sector, it is necessary for the builders to seriously consider the vital factors affecting the success of the project. In the meantime, it is necessary to allocate resources appropriately to the activities, and at the same time, there should be proper coordination between different team members. This coordination is considered necessary for teamwork and will have major and effective benefits in the process of implementing construction projects. Considering the special place of building information modeling (BIM) technology in the construction industry and its remarkable advantages in this industry, which facilitates the implementation process of the project, identifying interferences, and speeding up the implementation process, it is necessary to have different aspects of this technology. Furthermore, its role and position in the construction industry and its impact on building management should be investigated. Finally, the results obtained from this research show the positive impact of using building information modeling technology in risk management, safety, reducing energy consumption, and the possibility of intelligent project management.

## Introduction

In recent years, many countries spend a large proportion of their budget on the construction industry. According to the needs of every society and especially the young population in need of housing, different governments take into consideration the planning of providing the necessary housing and do their long-term planning accordingly. In this, having a written program and using the world's latest technologies has a special place [1]. Considering that the two essential factors of cost and time affect many construction projects, it is necessary to choose technology and implementation methods and techniques based on these two categories. At the same time, considering the prevailing competitive atmosphere in the construction industry and its economic aspect, both of which are growing and flourishing in parallel with each other. Therefore, the developed countries of the world are trying to gain a significant share of the market by introducing new technologies, taking over construction, and presenting themselves at the international level [2].

In the ever-evolving world of construction, implementing advanced technologies has become increasingly crucial to enhance project efficiency and productivity. Among these technologies, building information modeling is utilized in construction projects to systematically manage the various stages of work from design to operation with precision, elegance, and better management. By integrating various aspects like virtual design prototyping, 3D modeling, clash detection, and real-time collaboration, BIM enables project stakeholders to visualize and analyze the entire project life cycle before physical construction begins. One of the significant benefits of incorporating BIM into construction processes and performing teamwork will be identifying and eliminating possible errors in the design phase and not spreading them to the project's implementation phase, causing interruptions in work and exorbitant costs [3].

The construction industry plays a crucial role in shaping the infrastructure and development of nations. However, this industry is known as one of the most high-risk industries; thus, it is inevitable to pay sufficient attention to its safety category. Currently, due to the development of modern technologies, the traditional methods used in safety management have given way to the mentioned technologies. Therefore, using modern technologies such as BIM can effectively reduce human and financial losses while increasing productivity and improving construction quality. By providing a comprehensive virtual representation of the construction site, including architectural models, structural systems, and mechanical, electrical, and plumbing (MEP) designs, BIM allows stakeholders to visualize and address potential safety hazards even before physical work commences [4].

Adopting building information modeling technology and cost-saving measures, along with enhanced productivity, can lead to several advantages within the construction sector. Such benefits include reduced construction time, minimized material wastage, prevention of collisions and contradictions between structural and non-structural elements, as well as fostering effective coordination amongst structural engineers, architects, and facilities personnel. Furthermore, such practices facilitate teamwork, generate opportunities, and contribute to the creation of favorable conditions for modular construction. Additionally, they promote improved communication among personnel, the integration of data transfer systems, the mitigation of potential risks and hazards, and the possibility of constructing eco-friendly buildings that conserve energy [5].

One of the most crucial advantages of using this technology in the construction industry is that the project implementation process is done in a shorter period with higher quality and productivity. This method plays a significant role in improving the communication between the personnel working on the project and the exchange of information and data. At the same time, the variety of software available in BIM makes it possible for designers and builders to identify and fix existing interferences before entering the implementation phase, in the design phase, so that after the start of the construction process, significant savings in time and money can be made. Also, BIM provides the possibility of better and smarter selection of materials for employers. At the same time, by improving the planning phase, BIM will have a significant effect in preventing the wastage of materials. As a result, cost reduction, considered one of the most vital project management factors, is realized [6].

The objective of this paper is to examine the role and significance of incorporating building information modeling (BIM) across various sectors. Specifically, the investigation aims to illuminate how BIM can enhance smartness and safety within these sectors, as well as its potential to mitigate energy consumption. By comprehensively exploring these areas, this study contributes to the body of knowledge surrounding the utilization of BIM, thereby informing and advancing both academic and practical understanding of its multifaceted applications.

## 1. The effect of BIM from an organizational point of view

In general, the implementation of BIM technology can have a significant impact on organizational performance. For this purpose, individuals must first prepare their organizational goals and complete the ways to achieve these goals according to the progress of the work. Therefore, in all stages of the work, they must answer the question of how the use of BIM will help in achieving the project's goals [7].

From an organizational point of view, they should pay attention to the fact that BIM tools increase the technical capabilities of the project and help provide a suitable platform for providing special services to employers [8].

At the same time, the implementation of a comprehensive and intelligent project management system and BIM technologies within an organization can have clear or unclear benefits. In this regard, the organization is obliged to compile and specify the amount of information that is stored in the components of the model. The amount of information also depends on what they will be used for [9].

## 2. Main methods in using BIM

### 2.1. Central tank method

In this method, it is assumed that all project information is stored in a single database file. For example, all project scheduling and financial estimation information will be added to its 3D model information. This method could be more reasonable and practical. The type of information the designer requires differs from the contractor [10]. While the designer is involved in issues such as checking the energy consumption of the building, applying regulations, and designing the spaces, the contractor is interested in determining the work schedule and cost estimation; therefore, in order for the estimation work to begin, the designer's work must be finished. This is not practical, at least in the initial stages of work [11].

### 2.2. Extended reservoir methods

The extensive tank method is the method used by most designers and contractors. In this method, the BIM model has access to a set of separate databases created by independent programs. For example, all the information needed for the financial estimation of the project is in the corresponding independent program [12]. In order to do its work, this program needs two-way communication with the 3D BIM model to exchange the necessary information. This work is possible in the early stages of design; therefore, despite the use of independent data sources, due to a property called Interoperability, all the information of the different groups involved in the project is integrated. In this way, the design groups, including architecture, structure, electrical, and mechanical facilities, prepared their models separately in software such as Autodesk Revit. Finally, with the help of software such as Autodesk Navis Works, they were put together to obtain an integrated BIM model [13].

## 3. BIM deployment and implementation plan

The BIM implementation plan has relevant details that define how a project relying on BIM will be implemented, monitored, and organized. Therefore, the primary purpose of the BIM implementation plan, which will be displayed as BEP from now on, is to provide a general plan to ensure that all parties involved in their priorities and responsibilities in line with the implementation of BIM are aware. Therefore, BEP should be considered as a changeable document. At the same time, obtaining assurance from the project schedule's subjects is necessary [14].

It should be noted that some projects may have several BIM deployment plans so as to overlap different phases, including design, construction, and facility management (FM) [15].

## 4. Information sharing and cooperation of different team members

Usually, information sharing and collaboration between different members of the BIM team, such as the architect, design engineer, installation engineer, project manager, and other group members, are done through the Cloud Services system [16].

With the help of this cloud system, all the different parts of the BIM team are able to access the information, and any changes made by each one are automatically provided to the rest of the group members. At the same time, by using information modeling technology, companies and construction equipment manufacturers can quickly provide their information and modeled data to designers and experts [17].

## 5. Phased review of the building by BIM

In general, in the phased review of the building by BIM, there will be three main stages, which are:

- Step 1: before construction
- Step 2: under construction
- Step 3: after construction

It should be noted that all the above steps are interdependent, and there is no clear boundary between them. One can use these steps to come to life cycle design strategies that focus on minimizing the environmental impact of a building [18].

In the pre-construction phase, site selection, building design, and processes related to building design before establishment and operation include materials. Under the sustainable design strategy, the environmental consequences of structure design, building orientation, and its effects on the landscape and the materials used are examined [19].

The construction phase refers to a phase of the building's life cycle that is physically being built and operated. In sustainable design strategy, building construction and operation processes, in order to find ways to reduce harmful environmental impacts consumption of resources are examined; in addition, to the long-term effects of the building's surroundings on the health of its occupants is paying attention to the post-construction phase also begins when the useful life of the building is over. At this point, the case material used in the building will be used as resources for other buildings or will return to nature as waste [20].

Sustainable design strategy emphasizes reducing construction waste by recycling and reusing buildings and construction materials. It has special attention. In the division of other implementation stages of BIM, the work process is classified into five main subgroups, which are:

- First phase: pre-construction
- Second phase: construction
- Third phase: detailed design
- Fourth phase: preliminary design
- Fifth phase: planning

The effect of proper project management is in the design phase on the possibility of rework.

The cost spent on repair of damage and rework has become one of the significant issues in construction in most countries today. American Philip Krusi believes that more than 15 to 20 percent of the net profit of the construction is spent on reworking the construction so that the buildings have an acceptable quality [21].

In general, researchers believe that the successful integration of different building plans and the need to prepare BIM can be as simple as possible in designing, manufacturing, transporting, and delivering prefabricated parts, which has a direct impact on project management, among the other advantages of using this method, we can mention the reduction of direct costs, timely delivery of parts at the appointed time. Building Information Modeling (BIM) refers to a digital software collection that aims to facilitate project coordination and construction. BIM can also be considered as a process of documenting the stages of construction and project implementation. Before starting the construction process, the executed plan of the building should be prepared and built on the computer. In other words, BIM is a multi-dimensional model the fourth dimension of time, and the fifth dimension of cost. BIM is not only scalable in terms of project size and complexity but also in other ways. By using BIM tools, plans can be expanded directly in 3D space as a set of model elements. The advantages of using the BIM system are [22]:

### 5.1. Ensuring construction safety through BIM

Despite the economic aspect and the level of profitability of the construction industry in the world, this industry is always considered one of the riskiest industries. Although the number of injuries and deaths in the construction industry has decreased over the last two decades, this industry is still considered one of the most dangerous industries. According to the statistics announced by the Health and Safety Association (HSE) of England, about 141 injuries resulting in death occurred in 2012, 23 of which were related to the construction industry. About a third of the injuries that lead to death in the UK occur in the construction industry. According to global statistics, the construction industry faces more injuries yearly than other industries. Various studies indicate that 11-42% of accidents can be prevented in the design and construction phase by using safety requirements [23].

Therefore, in order to develop the safety assessment system, it is necessary to collect the required BIM information in the planning stage [24]. At the same time, several practical applications can be used at any stage of the construction process to improve safety and health. Therefore, in this situation, the mass of information can be used in the BIM environment. It is also possible to create links between different documents and improve the relationship between part of the information and the user guide and the possibility of filtering the available information.

### 5.2. The role of using BIM in reducing energy consumption and green buildings

Due to the positive effects of green buildings from residents' environmental, economic, health, and productivity perspectives, the use of this concept has been met with great acceptance worldwide. Therefore, according to the recent developments in the construction industry, valuable technologies have been presented to the engineering community so that designers, builders, investors, and other factors affecting the construction industry can benefit from countless benefits. By using them, they can improve the productivity and quality of construction. In general, green buildings are built to build environmentally friendly buildings, reduce environmental effects, and prevent energy wastage. In the last decade, attention to the human design of buildings has also been of particular importance, so the principle of human design is based on preserving natural conditions, urban design, and design based on human comfort [25].

In general, the environmental benefits of green buildings include improving and protecting biodiversity and ecosystems, improving water and air quality, reducing waste and wastage of materials, and, finally, preserving and restoring natural resources. At the same time, the integration of BIM in the construction industry has also benefited from an economic point of view, which includes reducing current and operational costs, improving the productivity of occupants and operators, enhancing asset value and profitability, and, finally, optimizing economic performance. The role of BIM in the meantime is to create an optimal opportunity to use design data for sustainable design as well as their analysis. Until 1333, when the concepts and technologies related to BIM were discussed for the first time, continuous efforts were made to investigate the possibility of using BIM in green buildings. Various functions of BIM, such as simulation of energy performance, analysis of lighting, and analysis of waste resulting from the destruction of buildings, were studied. Several BIM applications were proposed, and many aimed at integrating sustainability analysis and their usability at different stages of design, construction, and operation vectors were developed [26].

### 5.3. Using BIM in risk management

In risk management, identifying risks is usually considered the first step. At the same time, the classifications of risks that cause various risks in structures are usually considered as part of the risk identification process [27].

In general, there are two types of risks: internal and external. Regarding external risks [28] and their identification, the agents and personnel in the workshop should be referred to because the project monitoring team has more control over the mentioned risks. The range of internal risks and risks is more comprehensive than external ones, and their identification by the project monitoring team requires more accuracy and time. At the same time, the number of internal risks is higher than external risks, and their relationship is more complex. Internal risks are classified into two main categories, general and partial because some impact the entire or part of the project. In the first place, the analysis of project risks, including risk identification, risk classification, and determination of risk modification limits, must be consistent with the risk management strategy and, if necessary, be modified at different stages. Often, the process of correcting risks is of a preventive type, and it is necessary to be examined at different stages of project implementation and by different executive agents. According to ISO 31010:2009, risk management is a logical and systematic method that includes activities and processes related to identifying, evaluating, analyzing, recording, and eliminating risks [29].

Building information modeling reduces construction time, costs, and claims. Nevertheless, what happens if BIM is not carefully shared and explained to all construction members? In this section, there will be risks that are categorized as BIM risk [30].

## 5.4. Conflict between project stakeholders and contractors

When implementing BIM as an integral part of the project, the essential and critical issue is cooperation and access to the model by all people in the construction stages. Architectural engineers can use BIM to create premium construction projects and reduce costs and materials. At the same time, if the appropriate intellectual approach is not transferred to the beneficiaries when designing the building, then claims and problems will arise [31].

## 5.5. Unprincipled changes in the model

After the BIM is created and developed for project members, preventive measures should be taken to reduce the possibility of subsequent changes to the model. Currently, the risk of model change has become one of the fundamental challenges in large infrastructure projects. Therefore, consultants or contractors may change the model for their goals and interests. In recent years, the use of BIM in the construction industry has experienced significant growth. It has the potential to strengthen cooperation and communication, increase productivity and quality, and reduce project cost and delivery time. In order to overcome the obstacles in the traditional method of risk management, countless efforts have been made to use BIM technologies and BIM-related technologies that are used in project risk management. For example, the effectiveness of BIM as a systematic method for early identification and assessment of risks related to design and construction has been proven through 3D visualization, 4D planning, and 5D cost estimation. A computer system's spatial visualization and dynamic project modeling can significantly facilitate risk identification and initial communication. At the same time, it improves strategy and decision-making, time and cost management [32].

At the same time, in order to use the experiences related to previous projects in the risk management process, it is possible to prepare a comprehensive risk database and list all possible risks that affect the project implementation process. The aforementioned database facilitates the systematic understanding of all project risks. It helps the project executive team to link information and data related to project risks to real projects and increase decision-making power. Since the construction industry is dynamic with unexpected changes and risks, and it is necessary to update the input information daily, using a logical and principled method to classify information is inevitable. Currently, various tools are used to classify risks, including preparing a risk list (list of risks), risk matrix, risk maps, and risk breakdown structure or RBS. RBS is a flexible tool that can be easily updated. The main advantages of the risk breakdown structure include the following:

1- Increasing the overall understanding of risks and facilitating the communication of risks
2- Helping to identify the location of risks and create various solutions in order to eliminate risks
3- Providing an efficient structure to manage risk database and develop risk management software

So far, the primary approach used to develop RBS is extracting risk-related data from academic publications, project reports, and past project experience and classifying risk factors into a number of logical groups according to the sources of risk [33].

## Conclusions

Using the BIM method in construction projects has many advantages, including reducing execution time, more efficient workshop management, more coordination between different members, information exchange, and information sharing between different project design and implementation departments.

At the same time, due to the excessive details and executive details of construction projects such as bridges and buildings, which include connections, plates, reinforcements, and openings, in order to increase accuracy. It will be inevitable to use the BIM method in the work and to avoid ignoring some fine details. In addition to the complexity of some construction projects, these projects also have different construction and implementation processes, so due to the different materials used in different regions and the specific climatic conditions of each region, it is possible to use particular methods, techniques, and strategies to be used executively. The result of using information modeling technology in construction projects is better design, better construction, better operation, better execution quality, correct and logical project management, the possibility of making timely and appropriate decisions, avoiding rework, and execution errors.

Considering the comprehensive and distinguished capabilities of building information modeling technology (BIM) and the possibility of using it in different dimensions and directions of design, implementation, and safety, this technology will have significant advantages. One of the most essential applications of BIM is project risk management. In order to use the experiences related to previous projects in the risk management process, it is possible to prepare a comprehensive risk database and list all possible risks that affect the project implementation process. The results obtained from this research show the positive effect of using building information modeling technology in risk management, safety, reducing energy consumption, and the possibility of intelligent project management.